# BAYESIAN POISSON PROCESS PARTITION CALCULUS WITH AN APPLICATION TO BAYESIAN LÉVY MOVING AVERAGES

By Lancelot F. James[1]

*The Hong Kong University of Science and Technology*

This article develops, and describes how to use, results concerning disintegrations of Poisson random measures. These results are fashioned as simple tools that can be tailor-made to address inferential questions arising in a wide range of Bayesian nonparametric and spatial statistical models. The Poisson disintegration method is based on the formal statement of two results concerning a Laplace functional change of measure and a Poisson Palm/Fubini calculus in terms of random partitions of the integers $\{1,\ldots,n\}$. The techniques are analogous to, but much more general than, techniques for the Dirichlet process and weighted gamma process developed in [*Ann. Statist.* **12** (1984) 351–357] and [*Ann. Inst. Statist. Math.* **41** (1989) 227–245]. In order to illustrate the flexibility of the approach, large classes of random probability measures and random hazards or intensities which can be expressed as functionals of Poisson random measures are described. We describe a unified posterior analysis of classes of discrete random probability which identifies and exploits features common to all these models. The analysis circumvents many of the difficult issues involved in Bayesian nonparametric calculus, including a combinatorial component. This allows one to focus on the unique features of each process which are characterized via real valued functions $h$. The applicability of the technique is further illustrated by obtaining explicit posterior expressions for Lévy–Cox moving average processes within the general setting of multiplicative intensity models. In addition, novel computational procedures, similar to efficient procedures developed for the Dirichlet process, are briefly discussed for these models.

**1. Introduction.** Let $N$ denote a Poisson random measure on an arbitrary Polish space $\mathcal{W}$ characterized by its nonatomic sigma-finite mean

Received February 2003; revised September 2004.
[1]Supported in part by the RGC Grant HKUST-6159/02P and the DAG 01/02.BM43 of HKSAR.

intensity,

$$E[N(dw)] = \nu(dw).$$

That is to say, $N$ is a discrete random measure such that, for disjoint sets $A$ and $B$, $N(A)$ is independent of $N(B)$. Additionally, for each bounded set $B$, $N(B)$ is a Poisson random variable with finite mean $E[N(B)] = \nu(B)$. Following Daley and Vere-Jones [7], $N$ takes its values in the space of boundedly finite measures, say $\mathcal{M}$, equipped with an appropriate sigma-field $\mathcal{B}(\mathcal{M})$. Denote the law of $N$ as $\mathcal{P}(dN|\nu)$. Additionally, $BM(\mathcal{W})$ denotes the collection of Borel measurable functions of bounded support on $\mathcal{W}$. The class of nonnegative functions in $BM(\mathcal{W})$ is denoted as $BM_+(\mathcal{W})$. The law of $N$ is also uniquely characterized by its Laplace functional given by

$$(1) \quad \mathcal{L}_N(f|\nu) = \int_{\mathcal{M}} e^{-N(f)} \mathcal{P}(dN|\nu) = \exp\left(-\int_{\mathcal{W}} (1 - e^{-f(w)})\nu(dw)\right)$$

for each $f \in BM_+(\mathcal{W})$, where $N(f) = \int_{\mathcal{W}} f(w)N(dw)$. Note that the Laplace functional is well defined for all positive functions $f$. For additional information, see [22], Chapter 12. The Laplace functional, (1), will play a fundamental role in our analysis. An essential part of our presentation involves extensions of the following well-known disintegration for a joint measure of a point $W \in \mathcal{W}$ and $N$:

$$(2) \quad N(dW)\mathcal{P}(dN|\nu) = \mathcal{P}(dN|\nu, W)E[N(dW)] = \mathcal{P}(dN|\nu, W)\nu(dW),$$

where $E[N(dW)] = \int_{\mathcal{M}} N(dW)\mathcal{P}(dN|\nu)$ and $\mathcal{P}(dN|\nu, W)$ is a conditional distribution of $N$, given a point $W$, and coincides with the conditional law of the random measure

$$N + \delta_W,$$

where $N$ is $\mathcal{P}(dN|\nu)$ and $W$ is a fixed point. The result in (2) is equivalent to the Fubini theorem

$$(3) \quad \int_{\mathcal{M}} \left[\int_{\mathcal{W}} g(w, N)N(dw)\right]\mathcal{P}(dN|\nu) = \int_{\mathcal{W}} \left[\int_{\mathcal{M}} g(w, N)\mathcal{P}(dN|\nu, w)\right]\nu(dw),$$

for each measurable positive or integrable function $g$. Additionally, from the definition of $\mathcal{P}(dN|\nu, W)$, the following change of measure formula holds:

$$(4) \quad \begin{aligned} &\int_{\mathcal{W}} \left[\int_{\mathcal{M}} g(w, N)\mathcal{P}(dN|\nu, w)\right]\nu(dw) \\ &= \int_{\mathcal{W}} \int_{\mathcal{M}} g(w, N + \delta_w)\mathcal{P}(dN|\nu)\nu(dw). \end{aligned}$$

Within the framework of Palm calculus, the disintegration (2) is well known and may be found in [21] or [7], where $\mathcal{P}(dN|\nu, W)$ is an example of a Palm distribution. The representation (2) has been used extensively in a variety of important applications in probability; see, for instance, [35]. However, its use has been absent from the Bayesian nonparametrics literature. Note that since $N$ is not a random probability measure, $\mathcal{P}(dN|\nu, W)$ does not have the interpretation of a posterior distribution. However, the use of (2) is already enough to derive the posterior distribution of a variety of proper random probability measures when $n = 1$.

Random measures based on Poisson processes play an important role in spatial statistical analysis and Bayesian nonparametric statistics. In this work we will introduce a methodology we call a *Poisson process partition calculus* that provides a unified treatment of the otherwise formidable *posterior* analysis of such random measures. The idea appears in the unpublished manuscript of James [18], which discusses a variety of applications. Here, we will present a streamlined discussion which focuses specifically on methodology to deduce key properties of general classes of random probability measures and random intensities, analogous to those which make the Dirichlet process (see [13]) an attractive process for Bayesian non- and semiparametric analysis. The methodology consists of two components which will be described in more detail in Section 2. The first component is a Laplace functional/exponential change of measure formula for Poisson random measures, which can be seen as a form of functional exponential tilting or *Esscher* transform. The second is an extension of (2) in terms of partitions of the integers $\{1, \ldots, n\}$. One function of this extension is to allow one to bypass otherwise complex combinatorial arguments. In order to show explicitly the flexibility of our methods, we describe large classes of random probability measures in Section 1.1 which can be expressed as functionals of Poisson random measures. Additionally, in Section 1.1.1 we describe the structures of interest that are analogous to those for the Dirichlet process. Section 2 describes the elements of the Poisson process partition calculus. Section 3 discusses how to use the results in Section 2 to obtain the posterior analysis of the class of models described in Section 1.1. Section 4 presents a more explicit posterior analysis of a class of Lévy moving averages or hazard rates subject to a multiplicative intensity model. We also show, briefly, how this analysis leads to the development of computational procedures analogous to those used in Dirichlet process mixture models. Section 5 presents the formal details of the proof of Proposition 2.2.

1.1. *General discrete random probability measures and related concepts.* Let $h$ denote a strictly positive jointly measurable function on $\mathcal{W} \times \mathcal{M}$. One may define a general class of random probability measures, $P$, on $\mathcal{W}$ as follows:

(5) $$P(dw) = h(w, N)N(dw),$$

where $h$ is chosen such that $\int_{\mathcal{W}} h(w, N)N(dw) = 1$. The precise conditions on $h$ may also place restrictions on $\nu$. Note, however, that countable additivity of $P$ automatically follows from the additivity property of integrals with respect to $N$. Formally, we will consider random elements $W_1, \ldots, W_n | P$ which are i.i.d. with distribution $P$ and $P$ is defined in (5) with law, say $\mathcal{P}(dP|\nu)$, determined by a Poisson random measure $N$ with law $\mathcal{P}(dN|\nu)$. This gives a decomposition of the joint distribution of $(\mathbf{W}, P)$. We are interested in identifying explicitly the disintegration of this joint distribution in terms of the posterior distribution of $P|\mathbf{W}$, say $\pi(dP|\mathbf{W})$, and the exchangeable marginal distribution of $\mathbf{W}$ given by

$$(6) \qquad \left[\prod_{i=1}^{n} P(dW_i)\right]\mathcal{P}(dP|\nu) = \pi(dP|\mathbf{W}) \int_{\mathcal{M}} \left[\prod_{i=1}^{n} P(dW_i)\right]\mathcal{P}(dP|\nu).$$

In principle, the most difficult task is, of course, to obtain a clear expression for the posterior distribution $\pi(dP|\mathbf{W})$. This can be formidable for $n = 1$ and due to obvious nonconjugacy, and other issues to be discussed below, becomes more difficult for general $n$. However, explicit expressions for the marginal distribution and the posterior distribution are naturally linked. Hence, it is instructive to examine more closely the marginal distribution. By de Finetti's theorem, it is evident that the structure

$$(7) \qquad \mathbb{P}(d\mathbf{W}|\nu) := \int_{\mathcal{M}} \left[\prod_{i=1}^{n} P(dW_i)\right]\mathcal{P}(dP|\nu)$$

is exchangeable. It is a general analogue of the Blackwell and MacQueen [5] Pólya urn distribution. Moreover, this distribution is such that the random vector $\mathbf{W}$ possibly consists of ties and hence, the posterior distribution itself, $\pi(dP|\mathbf{W})$, also depends on ties. This suggests, as is natural for exchangeable structures (see the discussion in [25]), that the characterization of these quantities can involve a substantial combinatorial component. Here we discuss decompositions of (7) in terms of random partitions of the integers induced by these ties.

1.1.1. *Random partitions, EPPF, marginal distributions.* It is clear that there is a one-to-one correspondence between $\mathbf{W}$ and $(\mathbf{W}^*, \mathbf{p})$, where, using notation similar to Lo [30], $\mathbf{W}^* = (W_1^*, \ldots, W_{n(\mathbf{p})}^*)$ denotes the distinct values of $\mathbf{W}$ and $\mathbf{p} = \{C_1, \ldots, C_{n(\mathbf{p})}\}$ stands for a partition of $\{1, \ldots, n\}$ of size $n(\mathbf{p}) \leq n$ recording which observations are equal. The number of elements in the $j$th cell, $C_j := \{i : W_i = W_j^*\}$, of the partition is indicated by $e_j$, for $j = 1, \ldots, n(\mathbf{p})$, so that $\sum_{j=1}^{n(\mathbf{p})} e_j = n$. When it is necessary to emphasize a further dependence on $n$, we will also use the notation $e_{j,n} := e_j$. It follows that the marginal distribution of $\mathbf{W}$ can be expressed in terms of a conditional distribution of $\mathbf{W}|\mathbf{p}$, which is the same as a conditional distribution of the unique values $\mathbf{W}^*|\mathbf{p}$ and the marginal distribution of $\mathbf{p}$. The marginal distribution of $\mathbf{p}$, denoted as $\pi(\mathbf{p})$ or $p(e_1,\ldots,e_{n(\mathbf{p})})$, is an *exchangeable partition probability function* (EPPF), that is, a probability distribution on $\mathbf{p}$ which is exchangeable in its arguments and only depends on the size of each cell. The best known case of an EPPF is the variant of the Ewens sampling formula (ESF) associated with the Dirichlet process with total mass $\theta > 0$, given as

$$(8) \qquad p_\theta(e_1,\ldots,e_{n(\mathbf{p})}) = \frac{\theta^{n(\mathbf{p})}\Gamma(\theta)}{\Gamma(\theta+n)} \prod_{j=1}^{n(\mathbf{p})} \Gamma(e_j),$$

which was derived by Ewens [12] and Antoniak [3]. The EPPF can be interpreted as the distribution of the configuration of ties (clusters) among the $\mathbf{W}$. To understand this relationship further, note that, analogous to the case of the Dirichlet process, one can define the following probabilities relevant to (7). Suppose that $W_{n+1}$ is a newly observed variable. Then the probability that $W_{n+1}$ is distinct from the values $\mathbf{W}$, given $\mathbf{p}$, is

$$(9) \qquad \mathbb{P}(W_{n+1} \text{ is new }|\mathbf{p}) = q_{0,n} = \frac{p(e_1,\ldots,e_{n(\mathbf{p})},1)}{p(e_1,\ldots,e_{n(\mathbf{p})})},$$

and for $j = 1,\ldots,n(\mathbf{p})$, the probability that $W_{n+1} = W_j^*$, given $\mathbf{p}$, is

$$(10) \qquad \mathbb{P}(W_{n+1} = W_j^*|\mathbf{p}) = q_{j,n} = \frac{p(e_1,\ldots,e_j+1,\ldots,e_{n(\mathbf{p})})}{p(e_1,\ldots,e_{n(\mathbf{p})})}.$$

It is known that, for the case of (8), one has $q_{0,n} = \theta/(\theta+n)$ and $q_{j,n} = e_j/(\theta+n)$, which are the probabilities associated with the Chinese restaurant process (see [38], page 60) and the Blackwell–MacQueen prediction rule. In principle, one can use the probabilities in (9) and (10) to generate samples from $\mathbf{p}$, according to the EPPF, via a generalized Chinese restaurant process. See [15] for a discussion. However, we point out that, in general, unlike the case of the Dirichlet process, these probabilities are not the probabilities, say $\mathbb{P}(W_{n+1} = W_j^*|\mathbf{W})$ for $j=1,\ldots,n(\mathbf{p})$, which correspond to the appropriate prediction rule of $W_{n+1}|\mathbf{W}$. Rather, the following relationship holds: for $j=1,\ldots,n(\mathbf{p})$,

$$\mathbb{P}(W_{n+1} = W_j^*|\mathbf{p}) = \int_{\mathcal{W}^{n(\mathbf{p})}} \mathbb{P}(W_{n+1} = W_j^*|\mathbf{W})\pi(d\mathbf{W}^*|\mathbf{p}),$$

where $\pi(d\mathbf{W}^*|\mathbf{p})$ denotes the distribution of $\mathbf{W}|\mathbf{p}$ in terms of the unique values $\mathbf{W}^*$.

REMARK 1. The general EPPF concept is described in [36, 37, 38, 39], where a variety of applications are discussed. The notation $\sum_{\mathbf{p}}$ will be used to denote the sum over all possible partitions of the integers $\{1,\ldots,n\}$. A general discussion of the marginal structures $\mathbb{P}(d\mathbf{W}|\nu)$, such as that presented here, does not seem available. In the language of the theory of random measures, $\mathbb{P}(d\mathbf{W}|\nu)$ is also seen to be the $n$th moment measure of $P$. That is, one can use it to obtain the integer moments of $P$ and related quantities.

**2. Poisson process partition calculus.** So far we have pinpointed the type of structures we would like to obtain. However, what is missing is a systematic and easy mechanism to get at explicit expressions for these quantities. The idea of this paper is to focus on the utilization of (partition based) disintegration results related to the joint measure of $(\mathbf{W}, N)$ given by

$$(11) \quad \left[\prod_{i=1}^{n} N(dW_i)\right]\mathcal{P}(dN|\nu) = \mathcal{P}(dN|\nu, \mathbf{W}) \int_{\mathcal{M}} \left[\prod_{i=1}^{n} N(dW_i)\right]\mathcal{P}(dN|\nu).$$

The quantity (11) is not a proper distribution. However, it is this general form on the left-hand side which appears, explicitly or in augmented form, in all the models that will be discussed. The right-hand side, similar to that of (6), consists of a conditional distribution of $N|\mathbf{W}$, $\mathcal{P}(dN|\nu, \mathbf{W})$ and a sigma-finite marginal measure of $N$,

$$(12) \quad M(d\mathbf{W}|\nu) = \int_{\mathcal{M}} \left[\prod_{i=1}^{n} N(dW_i)\right]\mathcal{P}(dN|\nu),$$

which behaves in many respects like an exchangeable urn distribution and, importantly, can be expressed in terms of $(\mathbf{W}^*, \mathbf{p})$. These quantities are direct extensions of (2). The main purpose of this section is to describe two results concerning the Poisson process and the disintegration of (11) which are fashioned as simple tools that can be tailor-made to address inferential questions arising in a wide range of Bayesian nonparametric models.

2.1. *Basic tools.* First an exponential change of measure or disintegration formulae based on Laplace functionals is given below. This is a simple functional extension of an analogous result for Lévy processes on $\mathcal{R}$ or more generally, $\mathcal{R}^d$, which may be found in [27], Proposition 2.1.3. Such an operation is commonly called *exponential tilting*.

PROPOSITION 2.1. *For each $f \in BM_+(\mathcal{W})$ and each $g$ on $(\mathcal{M}, \mathcal{B}(\mathcal{M}))$,*

$$\int_{\mathcal{M}} g(N)e^{-N(f)}\mathcal{P}(dN|\nu) = \mathcal{L}_N(f|\nu) \int_{\mathcal{M}} g(N)\mathcal{P}(dN|e^{-f}\nu),$$

*where $\mathcal{P}(dN|e^{-f}\nu)$ is the law of a Poisson process with intensity $e^{-f(w)}\nu(dw)$. In other words, the following absolute continuity result holds: $e^{-N(f)}\mathcal{P}(dN|\nu) = \mathcal{L}_N(f|\nu)\mathcal{P}(dN|e^{-f}\nu)$. The result extends to any nonnegative measurable $f$ such that $\int_{\mathcal{W}}(1 - e^{-f(w)})\nu(dw) < \infty$.*

PROOF. By the unicity of Laplace functionals for random measures on $\mathcal{W}$, it suffices to check this result for the case $g(N) = e^{-N(h)}$ for $h \in BM_+(\mathcal{W})$. It follows that

$$\int_{\mathcal{M}} e^{-N(f+h)} \mathcal{P}(dN|\nu) = \mathcal{L}_N(f|\nu) \int_{\mathcal{M}} e^{-N(h)} \mathcal{P}_f(dN),$$

where, for the time being, $\mathcal{P}_f$ denotes some law on $N$. Simple algebra shows that

$$\int_{\mathcal{M}} e^{-N(h)} \mathcal{P}_f(dN) = \frac{\mathcal{L}_N(f+h|\nu)}{\mathcal{L}_N(f|\nu)}$$

and, hence, $\mathcal{P}_f(dN) = \mathcal{P}(dN|e^{-f}\nu)$. The extension holds by the same argument, since $\mathcal{L}_N(f|\nu) > 0$. □

Now, while indeed it is possible to use (2) repeatedly to analyze many of the models discussed in Section 1.1, such an analysis does not circumvent the need for what might be formidable combinatorial analysis. One may note, for instance, the nontrivial arguments used by Antoniak [3] to derive (8). With this in mind, the next result, in Proposition 2.2, gives a partition-based representation of (11) which serves to significantly simplify such derivations for more general models. We will delay a proof of Proposition 2.2 until Section 5. First, we formally identify the law $\mathcal{P}(dN|\nu, \mathbf{W})$ appearing in (11) as a conditional distribution of $N$, given the points $\mathbf{W}$, which is equivalent to the law of the random measure

$$(13) \qquad N_n^* = N + \sum_{j=1}^{n(\mathbf{p})} \delta_{W_j^*},$$

where $N$ is $\mathcal{P}(dN|\nu)$ independent of the points $\mathbf{W}$. Note, by definition, for any measurable function $g$ on $\mathcal{W} \times \mathcal{M}$, that $\mathcal{P}(dN|\nu, \mathbf{W})$ satisfies the following change of variable, as in the case for $n = 1$:

$$(14) \quad \int_{\mathcal{M}} g(\mathbf{W}, N) \mathcal{P}(dN|\nu, \mathbf{W}) = \int_{\mathcal{M}} g\bigg(\mathbf{W}, N + \sum_{j=1}^{n(\mathbf{p})} W_j^*\bigg) \mathcal{P}(dN|\nu).$$

Using (14), it follows that the conditional Laplace functional of $N$ with respect to $\mathcal{P}(dN|\nu, \mathbf{W})$ is

$$\int_{\mathcal{M}} e^{-N(f)} \mathcal{P}(dN|\nu, \mathbf{W}) = \bigg[\prod_{j=1}^{n(\mathbf{p})} e^{-f(W_j^*)}\bigg] \int_{\mathcal{M}} e^{-N(f)} \mathcal{P}(dN|\nu)$$

(15)

$$= \mathcal{L}_N(f|\nu) \prod_{j=1}^{n(\mathbf{p})} e^{-f(W_j^*)}.$$

We now present the formal partition based disintegration of (11).

PROPOSITION 2.2. *Suppose that $(\mathbf{W}, N)$ are measurable elements in the space $\mathcal{W}^n \times \mathcal{M}$ having the joint measure in (11), where $N$ is a Poisson random measure with sigma-finite nonatomic mean measure $\nu$. Then the following disintegration holds:*

$$\left[\prod_{i=1}^{n} N(dW_i)\right] \mathcal{P}(dN|\nu) = \mathcal{P}(dN|\nu, \mathbf{W}) \prod_{j=1}^{n(\mathbf{p})} \nu(dW_j^*),$$

*where $\mathcal{P}(dN|\nu, \mathbf{W})$ corresponds to the law of $N$ determined by (15) and is representable in distribution as (13). The moment measure is expressible via conditional moment measures as*

$$M(d\mathbf{W}|\nu) = \prod_{j=1}^{n(\mathbf{p})} \nu(dW_j^*) = \nu(dW_1) \prod_{i=2}^{n} \left[\nu(dW_i) + \sum_{j=1}^{n(\mathbf{p}_{i-1})} \delta_{W_j^*}(dW_i)\right],$$

*where $n(\mathbf{p}_{i-1})$ is the size of the partition $\mathbf{p}_{i-1}$ of $\{1, \ldots, i-1\}$ encoding the ties between $W_1, \ldots, W_{i-1}$.*

One can combine Proposition 2.1 and Proposition 2.2, yielding the following useful result which will be used in Section 4.

PROPOSITION 2.3. *Suppose that $(\mathbf{W}, N)$ are measurable elements in the space $\mathcal{W}^n \times \mathcal{M}$, where $N$ is a Poisson random measure with sigma-finite nonatomic mean measure $\nu$. Then for each nonnegative measurable $f$ such that $\int_{\mathcal{W}} (1 - e^{-f(w)}) \nu(dw) < \infty$, the following disintegration holds:*

$$\left[\prod_{i=1}^{n} N(dW_i)\right] e^{-N(f)} \mathcal{P}(dN|\nu)$$

$$= \mathcal{L}_N(f|\nu) \mathcal{P}(dN|e^{-f}\nu, \mathbf{W}) \prod_{j=1}^{n(\mathbf{p})} e^{-f(W_j^*)} \nu(dW_j^*).$$

$M(d\mathbf{W}|e^{-f}\nu) = \prod_{j=1}^{n(\mathbf{p})} e^{-f(W_j^*)} \nu(dW_j^*)$ *is the nth moment measure of a Poisson random measure with intensity $e^{-f(w)}\nu(dw)$.*

PROOF. The proof of this result follows by first applying Proposition 2.1 to get

$$\left[\prod_{i=1}^{n} N(dW_i)\right] e^{-N(f)} \mathcal{P}(dN|\nu) = \mathcal{L}_N(f|\nu) \left[\prod_{i=1}^{n} N(dW_i)\right] \mathcal{P}(dN|e^{-f}\nu).$$

Conclude the result by applying Proposition 2.2 with $e^{-f(w)}\nu(dw)$ in place of $\nu(dw)$. □

**3. Formal Bayesian methodology.** We now describe how to use the results in Section 2 to obtain desired results for models such as (6). First define

(16)
$$\psi_n(\mathbf{W}) = \int_{\mathcal{M}} \left[ \prod_{j=1}^{n(\mathbf{p})} [h(W_j^*, N)]^{e_j} \right] \mathcal{P}(dN|\nu, \mathbf{W})$$
$$= \int_{\mathcal{M}} \left[ \prod_{j=1}^{n(\mathbf{p})} [h(W_j^*, N_n^*)]^{e_j} \right] \mathcal{P}(dN|\nu).$$

Then an application of Proposition 2.2 yields the following result.

THEOREM 3.1. *Let $P$ denote a random probability defined as in (5), where $N$ is a Poisson random measure with intensity $\nu$. Let $\mathbf{W} = (W_1, W_2, \ldots, W_n)$ denote a vector of random elements on a Polish space $\mathcal{W}$ such that $W_1, \ldots, W_n | P$ are i.i.d. with distribution $P$. Then the following results hold:*

(i) *The posterior distribution of $N|\mathbf{W}$, $\pi(dN|\nu, \mathbf{W})$, corresponds to the conditional law of the random measure*

(17)
$$N_n^* = N + \sum_{j=1}^{n(\mathbf{p})} \delta_{W_j^*},$$

*where now $\pi^*(dN|\mathbf{W}) = [\psi_n(\mathbf{W})]^{-1} \mathcal{P}(dN|\nu) \prod_{j=1}^{n(\mathbf{p})} [h(W_j^*, N_n^*)]^{e_j}$ is the conditional law of $N$, in (17), given $\mathbf{W}$.*

(ii) *The posterior distribution of $P|\mathbf{W}$ is equivalent to the conditional distribution of the random probability measure*

$$P_n^*(dw) = h(w, N_n^*) N_n^*(dw) = h(w, N_n^*) N(dw) + \sum_{j=1}^{n(\mathbf{p})} h(W_j^*, N_n^*) \delta_{W_j^*}(dw),$$

*where the law of $N|\mathbf{W}$ is $\pi^*(dN|\mathbf{W})$*

(iii) *The joint exchangeable marginal distribution of $\mathbf{W}$ is given by $\mathbb{P}(d\mathbf{W}|\nu) = \psi_n(\mathbf{W}) \prod_{j=1}^{n(\mathbf{p})} \nu(dW_j^*)$. Additionally, the EPPF derived from the marginal distribution of $\mathbf{W}$ is expressible as*

(18)
$$p(e_1, \ldots, e_{n(\mathbf{p})}) = \int_{\mathcal{W}^{n(\mathbf{p})}} \psi_n(\mathbf{w}) \prod_{j=1}^{n(\mathbf{p})} \nu(dw_j^*).$$

PROOF. The key point to note is that, since $P$ is a functional of $N$, results for the joint distribution of $(\mathbf{W}, P)$ follow from the corresponding joint distribution of $(\mathbf{W}, N)$. From (6), the joint distribution of $(\mathbf{W}, N)$ is expressible as $[\prod_{i=1}^{n} h(W_i, N)][\prod_{i=1}^{n} N(dW_i)]\mathcal{P}(dN|\nu)$. Applying Proposition 2.2, along with the identity $\prod_{i=1}^{n} h(W_i, N) = \prod_{j=1}^{n(\mathbf{p})} [h(W_j^*, N)]^{e_j}$, the joint distribution of $(\mathbf{W}, N)$ can be expressed as

$$\tag{19} \left[\prod_{j=1}^{n(\mathbf{p})} [h(W_j^*, N)]^{e_j}\right] \mathcal{P}(dN|\nu, \mathbf{W}) \prod_{j=1}^{n(\mathbf{p})} \nu(dW_j^*).$$

One now only needs to apply simple Bayes rule to obtain an expression in terms of the posterior distribution of $N|\mathbf{W}$ and the marginal distribution of $\mathbf{W}$. Formally, to obtain the marginal distribution of $\mathbf{W}$, one integrates out $N$ in (19), yielding the form of $\mathbb{P}(d\mathbf{W}|\nu)$ in (iii). The expression in (18) is then evident. Now, since $\psi_n(\mathbf{W}) > 0$, it follows that the posterior distribution of $N|\mathbf{W}$ is $\pi(dN|\mathbf{W}) = [\psi_n(\mathbf{W})]^{-1}[\prod_{j=1}^{n(\mathbf{p})} [h(W_j^*, N)]^{e_j}]\mathcal{P}(dN|\nu, \mathbf{W})$. Statement (i) now follows by the change of measure formula (14). That is, the posterior Laplace function of $N|\mathbf{W}$ is

$$\int_{\mathcal{M}} e^{-N(f)} \pi(dN|\mathbf{W}) = \left[\int_{\mathcal{M}} e^{-N(f)} \pi^*(dN|\mathbf{W})\right] \prod_{j=1}^{n(\mathbf{p})} e^{-f(W_j^*)}.$$

Statement (ii) follows from the fact that $P(dw) = h(w, N)N(dw)$ and the representations of the posterior distribution of $N|\mathbf{W}$ in statement (i). □

REMARK 2. Statement (ii) describes the posterior distribution of $P|\mathbf{W}$ via the distribution of $P_n^*$ determined by $\pi^*(dN|\mathbf{W})$. As one application, the prediction rule of $W_{n+1}|\mathbf{W}$ can be readily computed as

$$\mathbb{P}(dW_{n+1}|\mathbf{W}) = \int_{\mathcal{M}} P_n^*(dW_{n+1})\pi^*(dN|\mathbf{W}).$$

3.1. *Discrete random probability measures defined by completely random measures.* The random probability measures defined in (5) are actually a bit different than the random probability measures commonly used in Bayesian nonparametrics. In particular, as we shall show, the class $P$ contains augmented forms of, say, the Dirichlet process or Doksum's [8] neutral to the right processes. In Bayesian nonparametrics many random probability measures are actually functionals of completely random measures (see [23, 26]), say, $\mu$ defined over a Polish space $\mathcal{Y}$. The class of completely random measures contains, for instance, the gamma process and the random hazard processes discussed in [14]. Completely random measures, ignoring fixed points of discontinuity, are representable in a distributional sense as functionals of Poisson random measures. We now describe this construction. Specify $\mathcal{W} = \mathcal{J} \times \mathcal{Y}$, where $\mathcal{J} = (0, \infty)$. Additionally, for points $w = (s, y)$, $N(ds, dy)$ denotes a Poisson random measure with mean intensity

$$E[N(ds, dy)] = \nu(ds, dy) = \rho(ds|y)\eta(dy).$$

Furthermore, it is assumed that $\rho$ and $\eta$ are selected such that, for each bounded set $B$ in $\mathcal{Y}$,

$$\tag{20} \int_B \int_{\mathcal{J}} \min(s,1)\rho(ds|y)\eta(dy) < \infty.$$

Now define a random measure $\mu$ on $\mathcal{Y}$ such that it may be represented in a distributional sense as

$$\tag{21} \mu(dy) = \int_{\mathcal{J}} sN(ds, dy).$$

Following Daley and Vere-Jones [7], the condition (20) guarantees that $\mu$ is in the space of boundedly finite measures $\mathcal{M}$ equipped with an appropriate sigma-field, $\mathcal{B}(\mathcal{M})$. If $\rho$ does not depend on $y$, then $\mu$ is said to be homogeneous. Furthermore, if $\mathcal{Y} = (0, \infty)$, then $\mu$ is sometimes called a subordinator. That is to say, a nonnegative Lévy process with stationary increments. Similar to the definition of $P$ in (5), one can define a general class of discrete random probability measures on $\mathcal{Y}$ as

$$\tag{22} P_\mu(dy) = q(y, \mu)\mu(dy) = q(y, \mu)\int_{\mathcal{J}} sN(ds, dy),$$

where $q$ is a strictly positive measurable function such that $P_\mu$ is a well-defined random probability measure. Note that the second representation in (22) reveals, via a natural augmentation, a class of random probability measures on $\mathcal{J} \times \mathcal{Y}$ defined as

$$\tag{23} \tilde{P}_\mu(ds, dy) = q(y, \mu)sN(ds, dy).$$

That is to say, $\tilde{P}_\mu(ds, dy)$ defined in (23) is a special case of (5) with the choice of $h(s, y, N) = sq(y, \mu)$.

Now set $W_i = (J_i, Y_i)$ for $i = 1, \ldots, n$ points in $\mathcal{J} \times \mathcal{Y}$ and denote the unique values as $W_j^* = (J_{j,n}, Y_j^*)$ for $j = 1, \ldots, n(\mathbf{p})$. Additionally, define a random measure

$$\mu_n^*(dy) = \int_0^\infty sN_n^*(dy) = \mu(dy) + \sum_{j=1}^{n(\mathbf{p})} J_{j,n}\delta_{Y_j^*}(dy).$$

Noting the form in (23), it follows that for $\mathbf{W} = (\mathbf{J}, \mathbf{Y})$,

$$\psi_n(\mathbf{J}, \mathbf{Y}) = \left[\prod_{j=1}^{n(\mathbf{p})} J_{j,n}^{e_j}\right]\phi_n(\mathbf{J}, \mathbf{Y}),$$

$$\text{where } \phi_n(\mathbf{J}, \mathbf{Y}) = \int_{\mathcal{M}} \left[\prod_{j=1}^{n(\mathbf{p})} [q(Y_j^*, \mu_n^*)]^{e_j}\right]\mathcal{P}(dN|\nu).$$

Additionally, let $\mathbf{s} = (s_1, \ldots, s_n)$ and $(s_{1,n}, \ldots, s_{n(\mathbf{p}),n})$ denote the arguments of $\mathbf{J} = (J_1, \ldots, J_n)$ and the collection $(J_{j,n})$, respectively. These facts lead to the following result.

THEOREM 3.2. *Let $P_\mu$ denote a random probability defined as in* (22), *where $N$ is a Poisson random measure on $\mathcal{W} = \mathcal{J} \times \mathcal{Y}$, with mean intensity $\nu(ds, dy) = \rho(ds|y)\eta(dy)$. Let $\mathbf{Y} = (Y_1, Y_2, \ldots, Y_n)$ denote a vector of random elements on $\mathcal{Y}$ such that $Y_1, \ldots, Y_n | P_\mu$ are i.i.d. with distribution $P_\mu$. Then the following results hold:*

(i) *The posterior distribution of $N|\mathbf{Y}$ corresponds to the conditional law of the random measure $N_n^* = N + \sum_{j=1}^{n(\mathbf{p})} \delta_{J_j, Y_j^*}$, where the conditional law of $N$ in this representation, given $\mathbf{J}, \mathbf{Y}$, is*

$$\pi^*(dN|\mathbf{J}, \mathbf{Y}) = [\phi_n(\mathbf{J}, \mathbf{Y})]^{-1} \left[ \prod_{j=1}^{n(\mathbf{p})} [q(Y_j^*, \mu_n^*)]^{e_j} \right] \mathcal{P}(dN|\nu).$$

*Additionally, the distribution of $\mathbf{J}|\mathbf{Y}$ is $\mathbb{P}(d\mathbf{J}|\mathbf{Y}, \nu) \propto \phi_n(\mathbf{J}, \mathbf{Y}) \prod_{j=1}^{n(\mathbf{p})} J_{j,n}^{e_j} \rho(dJ_{j,n}|Y_j^*)$. The law of $\mu_n^*(dy) = \int_0^\infty s N_n^*(ds, dy)$, given $\mathbf{Y}$, determined by the law of $N_n^*|\mathbf{Y}$, corresponds to the posterior distribution of $\mu|\mathbf{Y}$.*

(ii) *The posterior distribution of $P_\mu|\mathbf{Y}$ is equivalent to the conditional distribution, given $\mathbf{Y}$, of the random probability measure $P_{\mu_n^*}(dy) = q(y, \mu_n^*) \mu_n^*(dy)$.*

(iii) $\mathbb{P}(d\mathbf{Y}|\nu) = [\int_{\mathcal{J}^{n(\mathbf{p})}} \phi_n(\mathbf{s}, \mathbf{Y}) \prod_{j=1}^{n(\mathbf{p})} s_{j,n}^{e_j} \rho(ds_{j,n}|Y_j^*)] \prod_{j=1}^{n(\mathbf{p})} \eta(dY_j^*)$ *is the exchangeable marginal distribution of $\mathbf{Y}$. The EPPF derived from the marginal distribution of $\mathbf{Y}$ is expressible as*

$$(24) \quad p(e_1, \ldots, e_{n(\mathbf{p})}) = \int_{\mathcal{J}^{n(\mathbf{p})} \times \mathcal{Y}^{n(\mathbf{p})}} \phi_n(\mathbf{s}, \mathbf{y}) \prod_{j=1}^{n(\mathbf{p})} s_{j,n}^{e_j} \rho(ds_{j,n}|y_j^*) \eta(dy_j^*).$$

PROOF. First note the representation $\mu(dY_i) = \int_\mathcal{J} J_i N(dJ_i, dY_i)$ for $i = 1, \ldots, n$. Augmenting the joint distribution of $(\mathbf{Y}, P_\mu)$ by $\mathbf{J}$ yields the distribution of $(\mathbf{J}, \mathbf{Y}, \tilde{P}_\mu)$. Noting that $\mathbf{W} = (\mathbf{J}, \mathbf{Y})$, and using the identity $\prod_{i=1}^n J_i = \prod_{j=1}^{n(\mathbf{p})} J_{j,n}^{e_j}$, the posterior distribution of $N|\mathbf{J}, \mathbf{Y}$ and, hence, that of $\mu$ and $\tilde{P}_\mu$, follows directly from Theorem 3.1. Similarly, the joint distribution of $\mathbf{J}, \mathbf{Y}$ is given by statement (iii) of Theorem 3.1. This in turn yields the distributions of $\mathbf{J}|\mathbf{Y}$ and $\mathbf{Y}$. The distribution of $P_\mu$ follows from the fact that $P_\mu(dy) = \int_\mathcal{J} \tilde{P}_\mu(dy, ds)$. □

REMARK 3. The results in Theorems 3.1 and 3.2 serve the purpose of exploiting the common features of many random probability measures. This in turn allows one to avoid otherwise cumbersome intermediate arguments and focus on the unique features of each process. That is to say, similar to parametric Bayesian results obtained via classical Bayes rule, one will often require a finer analysis which now, given the results in Theorems 3.1 and 3.2, depends on exploiting the specific features of $h$ and $\nu$.

REMARK 4. If one sets $\rho(ds|y) := \rho(ds)$ such that $\int_0^\infty \rho(ds) = \infty$, and specifies $\eta(dy)$ to be a probability measure, then the choice of $h(s, y, N) = s/T$ for $T = \int_0^\infty \int_{\mathcal{Y}} sN(ds, dy) = \mu(\mathcal{Y})$ yields the homogeneous *Poisson–Kingman* random probability measures. This class has been discussed in varying generalities and contexts in, for instance, [18, 24, 35, 38, 39, 40]. The Dirichlet process with total mass $\theta$ arises by the choice of $\rho(ds) = \theta s^{-1} e^{-s} ds$. Using this choice, one can recover (8) from (24) or (18). More generally, using this choice of $h$, one obtains the EPPF given by Pitman [39].

REMARK 5. James [20] shows that Doksum's [8] neutral to the right processes can be obtained by the choice of $h(s, y, N) = se^{-Z(y-)}$, for $(s, y)$ in $[0, 1] \times (0, \infty)$, where $Z(y-) = \int_0^1 \int_0^\infty I_{\{x<y\}}[-\log(1-u)]N(du, dx)$, where now $\rho(ds|y)$ is a Lévy measure on $[0, 1]$ and $\eta$ is modeled as a cumulative hazard. The work of James [20] is an example of the type of refined analysis mentioned in Remark 3.

REMARK 6. One may define analogues of Dirichlet process mixture models (see [30]) by mixing $P$ or $P_\mu$ with a known density or probability mass function. The posterior analysis of such models follows as a simple consequence of Theorem 3.1 or Theorem 3.2 and Fubini's theorem. In particular, $\mathbb{P}(d\mathbf{W}|\nu)$ plays the role of a mixing measure, in analogy to the Blackwell–MacQueen distribution. A further generalization of these types of models is given in [33]. However, structurally such models are more closely related to models we will describe in the next section. That is to say, their analysis does not follow directly from Theorem 3.1 or Theorem 3.2.

**4. Multiplicative intensity models and Lévy–Cox moving averages.** Similar to Lo and Weng [32] (see also [10]), one can define random hazard rates or spatial intensities on a Polish space $\mathcal{X}$ as

$$(25) \qquad \lambda(x|\mu) = \int_{\mathcal{Y}} k(x|y)\mu(dy) = \int_{\mathcal{Y}} \int_0^\infty k(x|y)sN(ds, dy),$$

where $k(x|y)$ denotes a known positive measurable kernel on a Polish space $\mathcal{X} \times \mathcal{Y}$ assumed to be $\eta$-integrable over $\mathcal{Y}$. Additionally, $k$ is chosen such that, for a sigma-finite measure $\tau$ on $\mathcal{X}$ and each bounded set $B$, $\int_B k(x|y)\tau(dx) < \infty$ for each fixed $y$. Under this condition one may define a random cumulative intensity for each bounded set $B$ as

$$(26) \qquad \int_{\mathcal{Y}} \left[ \int_B k(x|y)\tau(dx) \right] \mu(dy).$$

The models (26) are also known as Lévy–Cox moving average models as discussed in [41, 42]. The models (25) can be used to model intensities of counting process models, or hazard rates of distribution functions. In particular, if $\mathcal{X} = (0, \infty)$, then one can define a random density $f$ as

$$(27) \qquad f(x|\lambda) = e^{-\Lambda(x)}\lambda(x) = S(x|\lambda)\lambda(x),$$

where $\Lambda(x) = \int_0^x \lambda(v)\,dv = \int_{\mathcal{Y}}[\int_0^x k(v|y)\,dv]\mu(dy)$ is a cumulative hazard and $S(x|\lambda) := e^{-\Lambda(x)}$ is the survival function denoting the probability that a random variable $X_1 \geq x$. We, of course, assume that $\Lambda(\infty) = \infty$. We will provide a detailed posterior analysis of the general class of Lévy moving averages assuming a multiplicative intensity likelihood, which we now describe. Suppose, as in [2], that, for each $i = 1, \ldots, m$, and fixed $\mu$, there is an independent counting process with mean intensity $\lambda(x)U_i(x)$, where $U_i(x)$ is a predictable process which is observable. We discuss some specific interpretations of this function below. Under this assumption the counting processes correspond to classes of multiplicative intensity models as discussed in [1]. Jacod [17] (see also [2, 32]) showed that the likelihood of such counting processes is absolutely continuous to the likelihood of Poisson process models. Here, for $n \leq m$, we work with the multiplicative intensity likelihood with a random intensity (25) which can be represented as

$$(28) \qquad L(\mathbf{X}|\mu) = e^{-\mu(g_m)}\prod_{i=1}^n \int_{\mathcal{Y}} k(X_i|Y_i)\mu(dY_i),$$

where $g_m(y) = \sum_{i=1}^m \int_{\mathcal{X}} U_i(x)k(x|y)\tau(dx)$ and, hence,

$$\mu(g_m) = \int_{\mathcal{X}}\left[\sum_{i=1}^m U_i(x)\right]\lambda(x)\tau(dx).$$

Note that throughout we assume that $k$ and $(U_i)$ are chosen such that $g_m$ is in $BM_+(\mathcal{Y})$. The model (28) suggests that there are $X_1, \ldots, X_n$ completely observed points and $m - n$ points, say $X_{n+1}, \ldots, X_m$, which are partially observed. Meanwhile, $\mathbf{Y} = (Y_1, \ldots, Y_n)$ can be viewed as missing data. The multiplicative intensity likelihood captures a large variety of models which appear in event history analysis. For example, if $\mathcal{X} = (0, \infty)$ and one sets $U_i(x) = I_{\{X_i \geq x\}}I_{\{x \in B_i\}}$ for a random set $B_i$ independent of $X_i$, then one can use this to model various censoring mechanisms. Specifically, setting $B_i = [0, D_i]$ for a random variable $D_i$ corresponds to a right censoring model. An extension to left truncation and right censored models is given by the choice $B_i = (V_i, D_i]$, where $V_i$ is a random variable almost surely less than $D_i$ (see [2], Section III.2). On the other hand, setting $\sum_{i=1}^m D_i(x) = 1$ leads to the likelihood of an inhomogeneous Poisson process with mean intensity $l(x)\tau(dx)$. Before proceeding to the posterior analysis, we first describe some more details about the special case of the class of random distributions defined by (27).

4.1. *Random hazard rates and densities.* Some specific examples of kernels $k$ used to define hazard rates $\lambda$ include the Dykstra and Laud [10] kernel, which corresponds to $k(x|y) = I_{\{y \leq x\}}$, where it follows that

$$
(29) \quad
\begin{aligned}
K(t|y) &:= \int_0^t k(x|y)\,dx = (t-y)I_{\{y \leq t\}} \quad \text{and} \\
\Lambda(t) &= \int_0^\infty (t-y)I_{\{y \leq t\}}\mu(dy)
\end{aligned}
$$

for $t \geq 0$. This choice of $k$ generates the family of nondecreasing hazard rates. Dragichi and Ramamoorthi [9] establish the consistency of this class of random hazard rates under wide choices of $\mu$. If one chooses an exponential kernel $k(x|y) = e^{-xy}$, then

$$K(t|y) = \int_0^t e^{-xy}\,dx = y^{-1}(1 - e^{-yt}) \quad \text{and} \quad \Lambda(t) = \int_0^\infty y^{-1}(1 - e^{-yt})\mu(dy).$$

As discussed in [32], this induces hazard rates which are completely monotone. See [34] for a variation of this model. If one is unsure of the shape of the hazard, then one can use any of the convolution kernels that one finds in classical kernel based density estimation, where, for $y = (m, \sigma) \in (-\infty, \infty) \times (0, \infty)$, a fairly simple choice is the rectangular kernel $k(x|m, \sigma) = I_{\{|x-m| \leq \sigma\}}$. See [16, 32, 41, 42] for various choices of $k$ on the real line and for spatial models. Notice that for a random variable $T$ the quantity $\lambda(t)$ represents the hazard rate of $T$ given $\mu$, that is,

$$\lambda(t)\,dt = \mathbb{P}(t \leq T < t + dt | T \geq t, \mu).$$

Note, however, that the quantity $E[\lambda(t)]$ does not have the interpretation as a prior specification for the hazard rate. For instance, in the case of the stable law of index $0 < \alpha < 1$, one has

$$E[\lambda(t)] = \int_{\mathcal{Y}} k(t|y)E[\mu(dy)] = \int_{\mathcal{Y}} k(t|y)\left[\int_0^\infty \frac{1}{\Gamma(1-\alpha)}s^{-\alpha}\,ds\right]\eta(dy) = \infty,$$

and we see that it is possible that $E[\lambda(t)] = \infty$ for all $t$. It follows that to appropriately evaluate the marginal hazard rate of $T$, one needs to first find the distribution of $\mu$ or $N$, given $T \geq t$. Setting $U_1(x) = I_{\{x<t\}}$, we have $g_1(y) = \int_0^t k(x|y)\,dx := K(t|y)$. Hence, setting $f_1(s,y) = g_1(y)s$, it follows that $S(t|\lambda) = e^{-N(f_1)}$ and an application of Proposition 2.1 gives

$$S(t)\mathcal{P}(dN|\nu) = \mathcal{P}(dN|e^{-f_1}\nu)E[S((t)|\lambda)],$$

where

$$E[S(t|\lambda)] = \mathcal{L}_N(f_1|\nu) = e^{-\int_{\mathcal{Y}}\int_0^\infty(1-e^{-sK(t|y)})\rho(ds|y)\eta(dy)}$$

denotes the marginal survival function of $T$. The quantity $\mathcal{P}(dN|e^{-f_1}\nu)$ denotes the law of a Poisson random measure with mean intensity $e^{-sK(t|y)}\rho(ds|y)\eta(dy)$ and represents the posterior distribution of $N|T \geq t$. The marginal hazard rate is obtained as

$$E[\lambda(t)|T \geq t] = \int_{\mathcal{M}} \lambda(t)\mathcal{P}(dN|e^{-f_1}\nu)$$
$$= \int_{\mathcal{Y}} k(t|y)\left[\int_0^{\infty} e^{-sK(t|y)}s\rho(ds|y)\right]\eta(dy).$$

In the stable case the marginal hazard rate becomes $\int_{\mathcal{Y}} k(t|y)[K(t|y)]^{\alpha-1}\eta(dy)$. Noting the specifications for the Dykstra and Laud kernel in (29), in the stable case with $\eta(dy) = dy$, the prior predictive hazard rate and survival function are

$$\lambda_{0,\alpha}(t|DL) = \alpha^{-1}t^{\alpha} \quad \text{and} \quad S_{0,\alpha}(t|DL) = e^{-(1/(\alpha(\alpha+1)))t^{\alpha+1}},$$

which corresponds to a Weibull distribution. We now show that a likelihood for this model based on right censored data is a special case of (28). Suppose that $T_1, \ldots, T_n|\mu$ are i.i.d. random variables with density $f(t|\lambda)$. Then their joint density can be expressed as $\prod_{i=1}^n f(T_i|\lambda) = \prod_{i=1}^n S(T_i)\lambda(T_i)$. If there are additionally $T_{n+1}, \ldots, T_m$ random times which are right censored by random times $D_{n+1}, \ldots, D_m$, that is, $T_l > D_l$ for $l = n+1, \ldots, m$, where we assume that the distribution of the censoring times does not depend on $\mu$, then the likelihood of $\mu$ based on $n$ completely observed times and $m - n$ right censored times takes the form

$$(30) \quad \left[\prod_{l=n+1}^m S(D_l|\lambda)\right]\prod_{i=1}^n S(T_i|\lambda)\lambda(T_i) = \left[\prod_{i=1}^m S(\min(T_i, D_i)|\lambda)\right]\prod_{i=1}^n \lambda(T_i),$$

where we set $D_i = \infty$ for $i = 1, \ldots, n$. Setting $U_i(x) = I_{\{T_i \geq x\}}I_{\{x \leq D_i\}} = I_{\{x \leq \min(T_i, D_i)\}}$ for $i = 1, \ldots, m$, one can write

$$\prod_{i=1}^m S(\min(T_i, D_i)|\lambda) = e^{-\mu(g_m)},$$

where, in this case, $g_m(y) = \sum_{i=1}^m \int_0^{\min(T_i, D_i)} k(x|y)\,dx$. Hence, it is not difficult to see that (30) is a special case of (28) with $\mu(g_m) = \sum_{i=1}^m \Lambda(\min(T_i, D_i))$.

4.2. *Posterior analysis of Lévy moving averages.* We now show how Proposition 2.3 is used to obtain the posterior distributional properties of the class of Lévy moving averages under the multiplicative intensity model. Here we actually focus on $\mu$. The approach used has similarities to that of Lo and Weng [32] in the case of weighted gamma processes. The analysis proceeds, as in the proof of Theorem 3.2, by introducing a suitable augmentation and then establishing the appropriate results for $N$. First, setting $f_{k,m}(s,y) = g_m(y)s$, it follows that $N(f_{k,m}) = \mu(g_m)$. We now provide some notation which will be used in the description of the posterior distribution. Throughout we assume, for integers $l, m$ and fixed $y$, the condition

$$(31) \qquad \kappa_l(e^{-f_{k,m}}\rho|y) = \int_0^\infty s^l e^{-g_m(y)s} \rho(ds|y) < \infty.$$

Define $C(\mathbf{X}) = \sum_{\mathbf{p}} \prod_{j=1}^{n(\mathbf{p})} \int_{\mathcal{Y}} [\prod_{i \in C_j} k(X_i|y)] \kappa_{e_j}(e^{-f_{k,m}}\rho|y)\eta(dy)$. Additionally, for $j = 1, \ldots, n(\mathbf{p})$, define distributions of the unique jumps $J_{j,n}$, each depending on a corresponding $Y_j^*$, as

$$(32) \qquad \mathbb{P}(J_{j,n} \in ds|Y_j^*) = \frac{s^{e_j} e^{-g_m(Y_j^*)s} \rho(ds|Y_j^*)}{\kappa_{e_j}(e^{-f_{k,m}^*}\rho|Y_j^*)}.$$

Using Proposition 2.3 and straightforward algebraic manipulations, that is, an appeal to Bayes rule, one arrives at the following description of the posterior distribution of $\mu$, given $\mathbf{X}$ and related quantities.

THEOREM 4.1. *Let $\mu(dy) = \int_0^\infty s N(ds, dy)$ denote a completely random measure on a Polish space $\mathcal{Y}$ with law determined by the law of the Poisson random measure $N$ with mean $\nu(ds, dy) = \rho(ds|y)\eta(dy)$ on $\mathcal{J} \times \mathcal{Y}$. Suppose that $\mathbf{X}|\mu$ has the multiplicative intensity likelihood specified in* (28). *Then the posterior distribution of $\mu|\mathbf{X}$ can be described in terms of the posterior distribution of $\mu|\mathbf{Y}, \mathbf{X}$ mixed over the posterior distribution of $\mathbf{Y}|\mathbf{X}$, which is described as follows:*

(i) *The posterior distribution of $N|\mathbf{Y}, \mathbf{X}$ is equivalent to the conditional law of the random measure $N_{n,m}^*(ds, dy) = N_{f_{k,m}}(ds, dy) + \sum_{j=1}^{n(\mathbf{p})} \delta_{J_{j,n}, Y_j^*}(ds, dy)$, where conditional on $(\mathbf{J}, \mathbf{Y}, \mathbf{X})$, $N_{f_{k,m}}$ is a Poisson random measure with intensity*

$$(33) \quad E[N_{f_{k,m}}(ds, dy)] = e^{-f_{k,m}(s,y)} \nu(ds, dy) = e^{-g_m(y)s} \rho(ds|y)\eta(dy),$$

*not depending on $(J_{j,n})$. Additionally, given $(\mathbf{Y}, \mathbf{X})$, the $(J_{j,n})$ are conditionally independent of $N_{f_{k,m}}$ and are mutually independent with each $J_{j,n}$ having the distribution depending on $Y_j^*$ specified in* (32).

(ii) *Statement* (i) *implies that $\mu|\mathbf{Y}, \mathbf{X}$ is equivalent to the conditional distribution, given $(\mathbf{Y}, \mathbf{X})$, of the random measure*

$$\mu_{n,m}^*(dy) = \int_0^\infty s N_{n,m}^*(ds, dy) = \mu_{g_m}(dy) + \sum_{j=1}^{n(\mathbf{p})} J_{j,n} \delta_{Y_j^*}(dy),$$

*where conditional on $\mathbf{Y}$ and $\mathbf{X}$, $\mu_{g_m}(dy) := \int_0^\infty s N_{f_{k,m}}(ds, dy)$ is a completely random measure with Lévy measure specified in* (33). *Additionally, the $(J_{j,n})$ are conditionally independent of $\mu_{g_m}$.*

(iii) *If $\lambda$ is a random hazard rate or intensity defined in* (25), *then its posterior distribution, given* $(\mathbf{Y}, \mathbf{X})$, *is equivalent to the conditional distribution of the random measure*

$$\lambda_{n,m}^*(x) = \int_{\mathcal{Y}} k(x|y)\mu_{g_m}(dy) + \sum_{j=1}^{n(\mathbf{p})} J_{j,n} k(x|Y_j^*).$$

(iv) *The conditional distribution of $\mathbf{Y}|\mathbf{X}$ can be expressed via the conditional distributions of $\mathbf{Y}|\mathbf{p}, \mathbf{X}$ and $\mathbf{p}|\mathbf{X}$ as follows: The distribution of $\mathbf{Y}|\mathbf{p}, \mathbf{X}$ is such that the unique values of $\mathbf{Y}$, $Y_1^*, \ldots, Y_{n(\mathbf{p})}^*$, are conditionally independent with distributions*

(34)
$$\mathbb{P}(dY_j^*|\mathbf{p}, \mathbf{X}) := \pi(dY_j^*|C_j) \propto \left[\prod_{i \in C_j} k(X_i|Y_j^*)\right] \kappa_{e_j}(e^{-f_{k,m}}\rho|Y_j^*)\eta(dY_j^*).$$

$\pi(\mathbf{p}|\mathbf{X}) = [C(\mathbf{X})]^{-1} \prod_{j=1}^{n(\mathbf{p})} \int_{\mathcal{Y}} [\prod_{i \in C_j} k(X_i|y)]\kappa_{e_j}(e^{-f_{k,m}}\rho|y)\eta(dy)$ *is the posterior distribution of $\mathbf{p}|\mathbf{X}$.*

PROOF. Similar to the proof of Theorem 3.2, we work with an (augmented) joint distribution of $(\mathbf{X}, N)$. Removing the integrals in $L(\mathbf{X}|\mu)$ and making appropriate substitutions, it follows that a distribution of $(\mathbf{J}, \mathbf{Y}, N, \mathbf{X})$ is proportional to

(35) $$e^{-N(f_{k,m})} \left[\prod_{i=1}^{n} k(X_i|Y_i)J_i\right]\left[\prod_{i=1}^{n} N(dJ_i, dY_i)\right]\mathcal{P}(dN|\nu).$$

Using the identity $\prod_{i=1}^{n} k(X_i|Y_i)J_i = \prod_{j=1}^{n(\mathbf{p})} [\prod_{i \in C_j} k(X_i|Y_j^*)]J_{j,n}^{e_{j,n}}$, combined with an application of Proposition 2.3 to (35), shows that the joint distribution of $(\mathbf{J}, \mathbf{Y}, N, \mathbf{X})$ is proportional to

(36)
$$\mathcal{L}_N(f_{k,m}|\nu)\mathcal{P}(dN|e^{-f_{k,m}}\nu, \mathbf{J}, \mathbf{Y}) \prod_{j=1}^{n(\mathbf{p})} \left[\prod_{i \in C_j} k(X_i|Y_j^*)\right]$$
$$\times J_{j,n}^{e_j} e^{-g_m(Y_j^*)J_{j,n}} \rho(dJ_{j,n}|Y_j^*)\eta(dY_j^*),$$

where $\mathcal{P}(dN|e^{-f_{k,m}}\nu, \mathbf{J}, \mathbf{Y})$ corresponds to the conditional law, given $(\mathbf{J}, \mathbf{Y}, \mathbf{X})$, of the random measure $N_{n,m}^*(ds, dy) = N_{f_{k,m}}(ds, dy) + \sum_{j=1}^{n(\mathbf{p})} \delta_{J_{j,n}, Y_j^*}(ds, dy)$ described in statement (i). The distribution of $\mathbf{J}|\mathbf{Y}, \mathbf{X}$ is then obtained by integrating out $N$ in (36) and applying Bayes rule, using the finiteness condition (31). A similar procedure yields the distributions of $\mathbf{Y}|\mathbf{X}$.  □

REMARK 7. Note that the law of $N_{f_{k,m}}$ is also determined by first applying Proposition 2.1 to (35) to obtain

$$e^{-\mu(g_m)}\mathcal{P}(dN|\nu) = \mathcal{P}(dN|e^{-f_{k,m}}\nu)\mathcal{L}_N(f_{k,m}|\nu).$$

See [20] for a similar type of calculation for spatial NTR processes. Notice also that, conditional on $(\mathbf{J}, \mathbf{Y}, \mathbf{X})$, the dependence of $N_{f_{k,m}}$ (and $N^*_{n,m}$) on $\mathbf{X}$ is only through the function $f_{k,m}$.

REMARK 8. The marginal distribution of $\mathbf{Y}|\mathbf{X}$ can also be written as

$$(37) \qquad \pi(d\mathbf{Y}|\mathbf{X}) = [C(\mathbf{X})]^{-1}\left[\prod_{i=1}^{n} k(X_i|Y_i)\right] M_\mu(d\mathbf{Y}|e^{-f_{k,m}}\nu),$$

where

$$(38) \quad \begin{aligned} M_\mu(d\mathbf{Y}|e^{-f_{k,m}}\nu) &= \prod_{j=1}^{n(\mathbf{p})} \kappa_{e_j}(e^{-f^*_{k,m}}\rho|Y^*_j)\eta(dY^*_j) \\ &= \int_\mathcal{M}\left[\prod_{i=1}^{n}\mu(dY_i)\right]\mathcal{P}(dN|e^{-f_{k,m}}\nu) \end{aligned}$$

assumes a role analogous to the Blackwell–MacQueen Pólya urn distribution in Dirichlet process mixture models. This viewpoint becomes important when designing computational procedures.

REMARK 9. James [19] gives results for semi-parametric weighted gamma process mixture models under more complex multiplicative intensity structures, that is, for cases where the kernel $k$ depends on a Euclidean parameter $\beta$, and $\beta$ has prior distribution $\pi(d\beta)$. A careful examination of that work, coupled with the results given here, provides an obvious way to obtain the corresponding result for the general processes, via a straightforward application of Bayes rule. A notable wrinkle is that the Laplace functionals will depend on $\beta$, and, hence, one does not have the cancellation of the semi-parametric version of $\mathcal{L}_N(f_{k,m}|\nu)$. A discussion of this is omitted for brevity. See [16] for further details in the case of the weighted gamma process.

4.3. *Posterior intensity rates and predictive hazards.* Similar to the case of Dirichlet process mixture models, many posterior quantities can be expressed in terms of functionals of the missing values $\mathbf{Y}$ or the partition $\mathbf{p}$. For example, the posterior intensity rate depends upon the posterior mean for $\mu$. From Theorem 4.1, it follows that the posterior mean of $\mu|\mathbf{X}, \mathbf{Y}$ is given by

$$(39) \quad E[\mu^*_{n,m}(dy)|\mathbf{X}, \mathbf{Y}] = \kappa_1(e^{-f_{k,m}}\rho|y)\eta(dy) + \sum_{j=1}^{n(\mathbf{p})} E[J_{j,n}|Y^*_j]\delta_{Y^*_j}(dy),$$

where

$$E[J_{j,n}|Y_j^*] = \frac{\kappa_{e_j+1}(e^{-f_{k,m}}\rho|Y_j^*)}{\kappa_{e_j}(e^{-f_{k,m}}\rho|Y_j^*)} = \frac{\int_0^\infty s^{e_j+1}e^{-g_m(Y_j^*)s}\rho(ds|Y_j^*)}{\int_0^\infty u^{e_j}e^{-g_m(Y_j^*)u}\rho(du|Y_j^*)}.$$

The quantity (39) is also the conditional moment measure of $\mu_{g_m}$, given $(\mathbf{Y}, \mathbf{X})$. Using these expressions, we obtain the following generalization of Lo and Weng ([32], Theorem 4.2).

COROLLARY 4.1. *Theorem* 4.1 *implies that the posterior expectation of the intensity* (25), *given* $\mathbf{X}$ *and* $\mathbf{Y}$, *is*

$$E[\lambda(x)|\mathbf{Y}, \mathbf{X}] = \int_{\mathcal{Y}} k(x|y)\kappa_1(e^{-f_{k,m}}\rho|y)\eta(dy) + \sum_{j=1}^{n(\mathbf{p})} k(x|Y_j^*)E(J_{j,n}|Y_j^*)$$

*and, hence, the posterior expectation given* $\mathbf{X}$ *is*

$$E[\lambda(x)|\mathbf{X}] = \sum_{\mathbf{p}} \left( \int_{\mathcal{Y}} k(x|y)\kappa_1(e^{-f_{k,m}}|y)\eta(dy) \right.$$
$$\left. + \sum_{j=1}^{n(\mathbf{p})} \int_{\mathcal{Y}} k(x|y)E[J_{j,n}|y]\pi(dy|C_j) \right) \pi(\mathbf{p}|\mathbf{X}).$$

Note, importantly, that a predictive hazard rate is defined as $E[\lambda(X_{n+1})|\mathbf{X}]$.

REMARK 10. Corollary 4.1 shows that the posterior mean for the intensity rate can be estimated from Monte Carlo draws involving only $\mathbf{p}$ and $\mathbf{Y}^*$. Thus, in problems where inference focuses on estimating the intensity, there is no need to draw values from the posterior of $\mu$. From a computational perspective this can greatly simplify algorithms.

4.4. *Monte Carlo procedures.* Ishwaran and James [16] show that efficient sampling schemes used to approximate the posterior distributional properties of Dirichlet process mixture models can be applied with some modification to sample the posterior distribution of mixtures of weighted gamma processes in the present setting. A key point was to note the similarities between the distribution of $\mathbf{Y}|\mathbf{X}$ for Dirichlet process models relative to the Blackwell–MacQueen urn and (37) in the case of the weighted gamma process. Lo and Weng [32] and Lo, Brunner and Chan [31] also exploited this idea. Here we note that the explicit expression of (38) and its description in Theorem 4.1, for general processes $\mu$, allows one to extend some of these procedures. First note that if one wants to sample $\mu|\mathbf{X}$, one can obtain a draw from $\mathbf{Y}|\mathbf{X}$ and then draw from the distribution of $\mu_{n,m}^*|\mathbf{X}, \mathbf{Y}$ described in (ii) of Theorem 4.1. Here we give some ideas on how to sample from $\mathbf{Y}|\mathbf{X}$, noting that steps such as draws from $\mu|\mathbf{Y}, \mathbf{X}$ are natural additions. For brevity, we only sketch out some details, focusing on identifying the relevant probabilities, as one can deduce the operational formalities either from [15, 16] or other relevant cited works. Note that (38) is the $n$th moment measure of a completely random measure with Lévy measure specified in (33). That is, (38) is the $n$th moment measure of $\mu_{g_m}$ described in (ii) of Theorem 4.1. It follows that (38) can also be represented via its conditional moment measures [see (39)] as

$$\kappa_1(e^{-f_{k,m}}\rho|Y_1)\eta(dY_1)\prod_{r=1}^{n-1}\left[\kappa_1(e^{-f_{k,m}}\rho|Y_{r+1})\eta(dY_{r+1}) + \sum_{j=1}^{n(\mathbf{p}_r)} \frac{\kappa_{1+e_{j,r}}(e^{-f_{k,m}}\rho|Y_j^*)}{\kappa_{e_{j,r}}(e^{-f_{k,m}}\rho|Y_j^*)}\delta_{Y_j^*}(dY_{r+1})\right],$$

where $\mathbf{p}_r = \{C_{1,r}, \ldots, C_{n(\mathbf{p}_r),r}\}$ is the partition of $\{1,\ldots,r\}$ encoding the ties in the first $r$ observations $\mathbf{Y}_r = (Y_1, \ldots, Y_r)$ and $e_{j,r}$ is the cardinality of $C_{j,r}$. In order to simulate $\mathbf{Y}$ from (37), one can construct an analogue of the Pólya urn Gibbs sampler of Escobar [11] or sequential importance sampler (SIS) of Liu [28] by working with a density constructed from $E[\lambda(x)|\mathbf{Y}, \mathbf{X}]$. These procedures are duals. We first describe the idea for the SIS procedure. This procedure samples $Y_1, \ldots, Y_n$ sequentially based on the conditional densities, for $r = 0, \ldots, n-1$,

$$(40) \qquad \mathbb{P}(Y_{r+1} \in dy|\mathbf{Y}_r, \mathbf{X}) = \frac{l_{0,r}}{c_r}\lambda_r(dy) + \sum_{j=1}^{n(\mathbf{p}_r)} \frac{l_{j,r}(Y_j^*)}{c_r}\delta_{Y_j^*}(dy),$$

where $\lambda_r(dy) \propto k(X_{r+1}|y)\kappa_1(e^{-f_{k,m}}\rho|y)\eta(dy)$ and

$$l_{0,r} = \int_{\mathcal{Y}} k(X_{r+1}|y)\kappa_1(e^{-f_{k,m}}\rho|y)\eta(dy)$$

and

$$l_{j,r}(Y_j^*) = k(X_{r+1}|Y_j^*)\frac{\kappa_{1+e_{j,r}}(e^{-f_{k,m}}\rho|Y_j^*)}{\kappa_{e_{j,r}}(e^{-f_{k,m}}\rho|Y_j^*)}.$$

Furthermore, $c_r = l_{0,r} + \sum_{j=1}^{n(\mathbf{p}_r)} l_{j,r}(Y_j^*)$. The importance weights for this scheme are $\prod_{r=1}^{n-1} c_r$. Now, for $r = 0, \ldots, n-1$, let $\mathbf{Y}_{-(r+1),n}$ denote the collection of $n-1$ random variables determined by removing $Y_{r+1}$ from $(Y_1, \ldots, Y_n)$. A general analogue of the Pólya urn Gibbs sampler for generating $Y_1, \ldots, Y_n$ is implemented by drawing values $Y_{r+1}$ from the probabilities $\mathbb{P}(Y_{r+1} \in dy|\mathbf{Y}_{-(r+1),n}, \mathbf{X})$ for $r = 0, \ldots, n-1$. These probabilities are defined analogously to (40), where $\mathbf{Y}_{-(r+1),n)}$ plays the role of $\mathbf{Y}_r$. See [16] for more details in the case of the weighted gamma process.

As in the case of Dirichlet process mixture models, the SIS and Gibbs sampling procedures described above are attractive as one does not need to perform complex integration. However, if integration is manageable, then, due to a Rao–Blackwellization argument, it is generally better to apply the following new variation of the general weighted Chinese restaurant algorithms discussed in [15, 31]. We will describe an SIS procedure which has a dual Gibbs sampling procedure analogous to the collapsed Gibbs samplers. The key to the procedure is to generate partitions $\mathbf{p}$ based on probabilities defined using the *predictive hazard rate*. That is, for $r = 0 \ldots, n - 1$, define

$$l(r) = l_{0,r} + \sum_{j=1}^{n(\mathbf{p})} l_{j,r},$$

where $l_{j,r} = \int_{\mathcal{Y}} l_{j,r}(y)\pi(dy|C_{j,r})$. The distribution $\pi(dy|C_{j,r})$ is the distribution for the $j$th unique value, given $C_{j,r}$, defined similarly to (34). The special case when $r = 0$ corresponds to

$$l(0) = \int_{\mathcal{Y}} k(X_1|y)\kappa_1(e^{-f_{k,m}}\rho|y)\eta(dy).$$

By Corollary 4.1 it follows that $l(r)$ is the predictive hazard rate given $X_1, \ldots, X_r$ and $\mathbf{p}_r$. From this, it is possible to define a sequential algorithm to generate an importance draw for $\mathbf{p}$ from the posterior. The method can be described in terms of $n$ customers who enter a restaurant sequentially, similar to the class of WCR algorithms. However, now the role played by the EPPF for random probability measures in such algorithms is replaced by cumulants, $\kappa$, arising from Lévy measures. The first customer is seated to a table with probability $l(0)/l(0) = 1$. Now at step $r + 1$, given a configuration $\mathbf{p}_r = \{C_{1,r}, \ldots, C_{n(\mathbf{p}_r),r}\}$ of the integers $\{1, \ldots, r\}$, one determines the partition $\mathbf{p}_{r+1}$ by noting whether a customer $r + 1$ sits at a new table or sits at one of the existing tables $C_{j,r}$ for $j = 1, \ldots, n(\mathbf{p}_r)$. The seating rule is defined as follows. To seat customer $r + 1$, sit him at an occupied table $C_{j,r}$ with probability $\Pr(\mathbf{p}_{r+1}|\mathbf{p}_r) = l(r)^{-1}l_{j,r}$, where $\mathbf{p}_{r+1} = \mathbf{p}_r \cup \{r + 1 \in C_{j,r}\}$ for $j = 1, \ldots, n(\mathbf{p}_r)$. Otherwise, customer $r + 1$ sits at a new unoccupied table $C_{n(\mathbf{p}_r+1)}$ with probability $\Pr(\mathbf{p}_{r+1}|\mathbf{p}_r) = l(r)^{-1}l_{0,r}$, where $\mathbf{p}_{r+1} = \mathbf{p}_r \cup C_{n(\mathbf{p}_r+1)}$. After $n$ customers are seated, the algorithm will yield a partition $\mathbf{p} = \{C_1, \ldots, C_{n(\mathbf{p})}\}$ of $\{1, \ldots, n\}$. By James ([18], Lemma 2.3), this partition has density $q(\mathbf{p})$ satisfying

$$L(\mathbf{p})q(\mathbf{p}) = \prod_{j=1}^{n(\mathbf{p})} \int_{\mathcal{Y}} \left[\prod_{i \in C_j} k(X_i|y)\right] \kappa_{e_{j,n}}(e^{-f_{k,m}}\rho|y)\eta(dy),$$

where $L(\mathbf{p}) = \prod_{r=1}^{n} l(r-1)$. In other words, for any integrable function $t(\mathbf{p})$,

$$\sum_{\mathbf{p}} t(\mathbf{p})\pi(\mathbf{p}|\mathbf{X}) = \frac{\sum_{\mathbf{p}} t(\mathbf{p})L(\mathbf{p})q(\mathbf{p})}{\sum_{\mathbf{p}} L(\mathbf{p})q(\mathbf{p})}.$$

Thus, $q(\mathbf{p})$ is an importance density for drawing posterior values $\mathbf{p}$ with importance values $L(\mathbf{p})$. This fact, combined with Theorem 4.1, now suggests a method for approximating posterior quantities from the multiplicative intensity model:

1. Draw $\mathbf{p} = \{C_1, \ldots, C_{n(\mathbf{p})}\}$ from $q(\mathbf{p})$. Condition on $\mathbf{p}$ and draw $Y_j^*$ independently from $\pi(dY_j^*|C_j)$ for $j = 1, \ldots, n(\mathbf{p})$.
2. Use the value for $\mathbf{Y}$ from step 1 to draw $\mu$ from $\mu|\mathbf{Y}, \mathbf{X}$. That is, draw $\mu$ from the random measure $\mu_{g_m} + \sum_{j=1}^{n(\mathbf{p})} J_{j,n} \delta_{Y_j^*}$.
3. To approximate the posterior law of a functional $g(\mu)$, run the previous steps $B$ times independently, obtaining values $\mu^{(b)}$ with importance weights $L(\mathbf{p}^{(b)})$, for $b = 1, \ldots, B$. Approximate the law $\mathcal{P}\{g(\mu) \in \cdot|\mathbf{X}\}$ with

$$\frac{\sum_{b=1}^{B} I\{g(\mu^{(b)}) \in \cdot\} L(\mathbf{p}^{(b)})}{\sum_{b=1}^{B} L(\mathbf{p}^{(b)})}.$$

To approximate functions of the form $t(\mathbf{p})$, for instance, $E[\lambda(t)|\mathbf{X}]$, then step 2 can be eliminated and estimation is based on

$$\frac{\sum_{b=1}^{B} t(\mathbf{p}^{(b)}) L(\mathbf{p}^{(b)})}{\sum_{b=1}^{B} L(\mathbf{p}^{(b)})}.$$

REMARK 11. Note that the main difficulty in step 2 is to approximate a draw from $\mu_{g_m}$. There are several methods discussed in the literature. See, for instance, [4, 6, 42] for some possible ideas and further references in the general setting.

We next present some explicit examples of the posterior distribution of $\mu$ based on the results in Theorem 4.1.

4.4.1. *Generalized gamma process.* Brix [6] proposes an interesting class of measures by specifying $\mu$ to be a generalized gamma random measure. Using the description of Brix [6], these are $\mu$ processes with Lévy measure

$$\rho_{\alpha,b}(ds) = \frac{1}{\Gamma(1-\alpha)} s^{-\alpha-1} e^{-bs} ds.$$

The values for $\alpha$ and $b$ are restricted to satisfy $0 < \alpha < 1$ and $0 \leq b < \infty$ or $-\infty < \alpha \leq 0$ and $0 < b < \infty$. Different choices for $\alpha$ and $b$ in $\rho_{\alpha,b}$ yield various subordinators. These include the stable subordinator when $b = 0$, the gamma process subordinator when $\alpha = 0$ and the inverse-Gaussian subordinator when $\alpha = 1/2$ and $b > 0$. When $\alpha < 0$, this results in a class of gamma compound Poisson processes. Nieto–Barajas and Walker [34] provide analysis for a random distribution function on $(0, \infty)$, as in (27), where $k$ is an exponential kernel and where $\mu$ is modeled as a weighted version of a gamma compound Poisson process. This turns out to be an inhomogeneous variation of a subclass of the models of Brix [6] with $\alpha = -1$ and $b = b(y)$ in $BM_+(\mathcal{Y})$. The weighted gamma process considered in [32] corresponds to the choice of $\alpha = 0$ and $b = b(y)$.

The posterior distribution of $\mu$, given $(\mathbf{X}, \mathbf{Y})$, is equivalent to the conditional distribution of the random measure $\mu_{g_m} + \sum_{j=1}^{n(\mathbf{p})} (b + g_m(Y_j^*))^{-1} G_{j,n} \delta_{Y_j^*}$, where $\mu_{g_m}$ is an inhomogeneous generalized gamma process with intensity

$$\frac{1}{\Gamma(1-\alpha)} e^{-(g_m(y)+b)s} s^{-\alpha-1} \, ds \, \eta(dy),$$

and $(G_{j,n})$ are independent gamma random variables with shape $e_{j,n} - \alpha$ and unit scale. It follows that the conditional moment measure is

$$E[\mu_{n,m}^*(dy)|\mathbf{X}, \mathbf{Y}] = (b + g_m(y))^{\alpha-1} \eta(dy)$$
$$+ \sum_{j=1}^{n(\mathbf{p})} (b + g_m(Y_j^*))^{-1} (e_{j,n} - \alpha) \delta_{Y_j^*}(dy).$$

The joint moment measure of $\mathbf{Y}$ can be expressed as

$$M_\mu(d\mathbf{Y}|\rho_{\alpha, b+g_m} \eta) = \left[ \prod_{j=1}^{n(\mathbf{p})} \frac{\Gamma(e_{j,n} - \alpha)}{\Gamma(1-\alpha)} \right] \prod_{j=1}^{n(\mathbf{p})} (b + g_m(Y_j^*))^{-(e_{j,n}-\alpha)} \eta(dY_j^*),$$

which, if $b = b(y)$, generalizes an expression for the weighted gamma process; see [19, 32]. Note that, for $r = 0, 1, \ldots, n-1$,

$$l_{0,r} = \int_\mathcal{Y} k(X_{r+1}|y)(b + g_m(y))^{\alpha-1} \eta(dy)$$

and

$$l_{j,r}(Y_j^*) = \frac{k(X_{r+1}|Y_j^*)}{(b + g_m(Y_j^*))} (e_{j,r} - \alpha).$$

4.4.2. *Smoothed spatial beta process.* Given the conjugacy properties of the beta process when used as a cumulative hazard prior in [14] under right censoring, it is natural to think of a smooth version of this process to model hazard rates. This is in analogy to smoothing the Nelson–Aalen estimator. Here we allow an extension to $Y = (Y_1, Y_2) \in (0, \infty) \times \mathcal{Y}_2$ by specifying

$$\rho(ds|y_1) = c(y_1) s^{-1} (1-s)^{c(y_1)-1} \, ds$$

and writing $\eta(dy_1, dy_2)$, where $c$ is some positive function. Note, however, that the posterior behavior is quite different in this context than in [14]. The measure $\mu_{g_m}$ corresponds to a completely random measure with Lévy measure $c(y_1)e^{-g_m(y)s}s^{-1}(1-s)^{c(y_1)-1}\,ds\,\eta(dy_1, dy_2)$ and, hence, is not a beta process. Additionally, the distribution of $J_{j,n}$ is

$$\mathbb{P}(J_{j,n}|Y_j^*) = \frac{e^{-g_m(Y_j^*)s}s^{e_{j,n}-1}(1-s)^{c(Y_{1,j}^*)-1}\,ds}{\int_0^1 e^{-g_m(Y_j^*)u}u^{e_{j,n}-1}(1-u)^{c(Y_{1,j}^*)-1}\,du},$$

where the normalizing constant depends on the Laplace transform of a beta random variable evaluated at $g_m(Y_j^*)$. In other words, it is related to the *confluent hypergeometric function*

$$\begin{aligned}&{}_1F_1(e_{j,n}, c(Y_{1,j}^*) + e_{j,n}, -g_m(Y_j^*))\\&= \frac{\Gamma(c(Y_{1,j}^*) + e_{j,n})}{\Gamma(c(Y_{1,j}^*))\Gamma(e_{j,n})}\int_0^1 e^{-g_m(Y_j^*)u}u^{e_{j,n}-1}(1-u)^{c(Y_{1,j}^*)-1}\,du.\end{aligned}$$

For some simplification, hereafter we set $c$ equal to the constant $\theta$. Then it follows that one can write $E[\mu_{n,m}^*(dy)|\mathbf{Y}, \mathbf{X}]$ as

$$\theta\left[\int_0^1 e^{-g_m(y)s}(1-s)^{\theta-1}\,ds\right]\eta(dy)$$
$$+ \sum_{j=1}^{n(\mathbf{p})} \frac{e_{j,n}}{e_{j,n}+\theta} \frac{{}_1F_1(e_{j,n}+1, \theta+e_{j,n}+1, -g_m(Y_j^*))}{{}_1F_1(e_j, \theta+e_j, -g_m(Y_j^*))}\delta_{Y_j^*}(dy),$$

and the joint marginal measure $M_\mu(d\mathbf{Y}|e^{-f_{k,m}}\nu)$ is

$$\left[\prod_{j=1}^{n(\mathbf{p})} \frac{\Gamma(e_{j,n})\Gamma(\theta)}{\Gamma(e_{j,n}+\theta)}\right]\prod_{j=1}^{n(\mathbf{p})} {}_1F_1(e_{j,n}, \theta+e_{j,n}, -g_m(Y_j^*))\eta(dY_j^*).$$

**5. Proof of Proposition 2.2.** In this section we present two results which when combined lead to a proof of Proposition 2.2.

PROPOSITION 5.1. *Suppose the $(\mathbf{W}, N)$ are measurable elements in the space $\mathcal{W}^n \times \mathcal{M}$ having the joint measure in* (11), *where $N$ is a Poisson random measure with sigma-finite nonatomic mean measure $\nu$. Then the following disintegration holds:*

$$\begin{aligned}(41)\quad &\left[\prod_{i=1}^n N(dW_i)\right]\mathcal{P}(dN|\nu)\\&= \mathcal{P}(dN|\nu, \mathbf{W})\nu(dW_1)\prod_{i=2}^n\left[\nu(dW_i) + \sum_{j=1}^{n(\mathbf{p}_{i-1})}\delta_{W_j^*}(dW_i)\right],\end{aligned}$$

where $\mathcal{P}(dN|\nu,\mathbf{W})$ corresponds to the law of $N$ determined by (15) and is representable in distribution as (13). The statement implies that $M(d\mathbf{W}|\nu) = \nu(dW_1)\prod_{i=2}^{n}[\nu(dW_i) + \sum_{j=1}^{n(\mathbf{p}_{i-1})}\delta_{W_j^*}(dW_i)]$.

PROOF. First note the equivalence for $M(d\mathbf{W}|\nu)$ follows by integrating out $N$ in (41). The result proceeds by induction. The case for $n=1$, (2), is true. Now assuming that the result is true for $n=r$, it follows that

$$\left[\prod_{i=1}^{r+1} N(dW_i)\right]\mathcal{P}(dN|\nu) = N(dW_{r+1})\mathcal{P}(dN|\nu,\mathbf{W}_r)M(d\mathbf{W}_r|\nu),$$

which implies the form of $M(d\mathbf{W}_{r+1}|\nu)$, and, hence, it remains to show that

$$N(dW_{r+1})\mathcal{P}(dN|\nu,\mathbf{W}_r) = \mathcal{P}(dN|\nu,\mathbf{W}_{r+1})\left[\nu(dW_{r+1}) + \sum_{j=1}^{n(\mathbf{p}_r)}\delta_{W_j^*}(dW_{r+1})\right].$$

First, for functions $s$ and $f$ in $BM_+(\mathcal{W})$, note that, by a change of measure,

$$\begin{aligned}(42)\quad &\int_{\mathcal{M}}\int_{\mathcal{W}} s(w)e^{-N(f)}N(dw)\mathcal{P}(dN|\nu,\mathbf{W}_r) \\ &= \left[\prod_{j=1}^{n(\mathbf{p}_r)} e^{-f(W_j^*)}\right]\int_{\mathcal{M}} g(N_n^*)e^{-N(f)}\mathcal{P}(dN|\nu),\end{aligned}$$

where $g(N_n^*) = \int_{\mathcal{W}} s(w)N_n^*(dw) = \int_{\mathcal{W}} s(w)N(dw) + \sum_{j=1}^{n(\mathbf{p}_r)} s(W_j^*)$. Applying Proposition 2.1 to the right-hand side of (42) shows that the expressions in (42) are equal to

$$\mathcal{L}_N(f|\nu)\left[\prod_{j=1}^{n(\mathbf{p}_r)} e^{-f(W_j^*)}\right]\left[\int_{\mathcal{M}}\int_{\mathcal{W}} s(w)N(w)\mathcal{P}(dN|e^{-f}\nu) + \sum_{j=1}^{n(\mathbf{p}_r)} s(W_j^*)\right].$$

It follows that the conditional Laplace functional of $N$, given $\mathbf{W}_{r+1} := (\mathbf{W}_r, W_{r+1})$, relative to $M(d\mathbf{W}_{r+1}|\nu)$, is determined by the expression

$$\mathcal{L}_N(f|\nu)\left[\prod_{j=1}^{n(\mathbf{p}_r)} e^{-f(W_j^*)}\right]\left[\int_{\mathcal{W}} s(W_{r+1})e^{-f(W_{r+1})}\nu(dW_{r+1}) + \sum_{j=1}^{n(\mathbf{p}_r)} s(W_j^*)\right].$$

Now define a function $t(W_{r+1})$ to be $e^{-f(W_{r+1})}$ if $W_{r+1}$ is not equal to any of the $\{W_1^*,\ldots,W_{n(\mathbf{p}_r)}^*\}$ and is set to be one otherwise. Then, since $\nu$ is nonatomic, it follows that

$$\begin{aligned}&\int_{\mathcal{W}} s(W_{r+1})t(W_{r+1})\left[\nu(dW_{r+1}) + \sum_{j=1}^{n(\mathbf{p}_r)}\delta_{W_j^*}(dW_{r+1})\right] \\ &= \int_{\mathcal{W}} s(W_{r+1})e^{-f(W_{r+1})}\nu(dW_{r+1}) + \sum_{j=1}^{n(\mathbf{p}_r)} s(W_j^*).\end{aligned}$$

Hence, the conditional Laplace functional of $N$, given $\mathbf{W}_{r+1}$, with respect to $M(d\mathbf{W}_{r+1}|\nu)$ is

$$(43) \quad \mathcal{L}_N(f|\nu)\left[\prod_{j=1}^{n(\mathbf{p}_r)} e^{-f(W_j^*)}\right] t(W_{r+1}) = \mathcal{L}_N(f|\nu)\left[\prod_{j=1}^{n(\mathbf{p}_{r+1})} e^{-f(W_j^*)}\right],$$

as desired. □

The next result, which builds on Proposition 5.1, establishes the partition representation of $M(d\mathbf{W}|\nu)$.

PROPOSITION 5.2. *For $i = 1,\ldots,n$, let $g_i$ be nonnegative functions in $BM(\mathcal{W})$. Then*

$$(44) \quad \int_\mathcal{M}\left[\prod_{i=1}^n \int_\mathcal{W} g_i(w_i) N(dw_i)\right]\mathcal{P}(dN|\nu) = \sum_\mathbf{p} \prod_{j=1}^{n(\mathbf{p})} \int_\mathcal{W}\left[\prod_{i \in C_j} g_i(w_j^*)\right]\nu(dw_j^*).$$

*Equivalently, $M(d\mathbf{W}|\nu) = \prod_{j=1}^{n(\mathbf{p})} \nu(dW_j^*)$.*

PROOF. The proof of (44) proceeds by induction. Case $n = 1$ is obvious. Now suppose it is true for $n = r$. Let $\mathbf{p}_{r+1}$ denote a partition of $\{1,\ldots,r+1\}$, and define, for each $r > 0$,

$$\phi_g(\mathbf{p}_r) = \prod_{j=1}^{n(\mathbf{p}_r)} \int_\mathcal{W}\left[\prod_{i \in C_{j,r}} g_i(w_j^*)\right]\nu(dw_j^*).$$

It follows that $\phi_g(\mathbf{p}_{r+1})$ is $\phi_g(\mathbf{p}_r) \int_\mathcal{W} g_{r+1}(v)\nu(dv)$ if $n(\mathbf{p}_{r+1}) = n(\mathbf{p}_r) + 1$. Otherwise, if the index $r+1$ is in an existing cell/table $C_{i,r}$, then it is equivalent to $\phi_g(\mathbf{p}_r) \int_\mathcal{W} g_{r+1}(v)\pi_g(dv|C_{i,r})$, where

$$\pi_g(dv|C_{i,r}) = \frac{[\prod_{l \in C_{i,r}} g_l(v)]\nu(dv)}{\int_\mathcal{W}[\prod_{l \in C_{i,r}} g_l(v)]\nu(dv)}$$

for $i = 1,\ldots,n(\mathbf{p}_r)$. Note that this implies that

$$\sum_{\mathbf{p}_{r+1}} \phi_g(\mathbf{p}_{r+1}) = \sum_{\mathbf{p}_r} \phi_g(\mathbf{p}_r)\left[\int_\mathcal{W} g_{r+1}(v)\nu(dv) + \sum_{i=1}^{n(\mathbf{p}_r)} \int_\mathcal{W} g_{r+1}(v)\pi_g(dv|C_{i,r})\right].$$

Now, by (simple algebra) and the induction hypothesis on $r$, it follows that

$$\sum_{\mathbf{p}_{r+1}} \phi_g(\mathbf{p}_{r+1})$$

$$= \int_{\mathcal{W}^n}\left[\int_\mathcal{W} g_{r+1}(v)\nu(dv) + \sum_{j=1}^{n(\mathbf{p}_r)} g_{r+1}(W_j^*)\right]\left[\prod_{i=1}^r g_i(W_i)\right] M(d\mathbf{W}_r|\nu).$$

Now, utilizing the fact that $M(d\mathbf{W}_{r+1}|\nu) = [\nu(dW_{r+1}) + \sum_{j=1}^{n(\mathbf{P}_r)} \delta_{W_j^*}(dW_{r+1})] \times M(d\mathbf{W}_r|\nu)$ concludes the proof. Note this last statement relies on the result in Proposition 5.1. $\square$

REMARK 12. The proof of Proposition 5.2 follows closely an unpublished proof by Albert Lo for the case of gamma processes. That is, it is an alternative proof for Lemma 2 in [30] which yields the appropriate partition representation for integrals with respect to a Blackwell–MacQueen urn distribution derived from a Dirichlet process. The style of proof exploits properties of partitions similar to those stated in [36], Proposition 10. Details in the proof of Proposition 5.2 translate into justifications for generalizations of weighted Chinese restaurant algorithms.

**Acknowledgments.** I would like to thank Jim Pitman for convincing me to look at the general Poisson framework, and generously sharing his expertise with me. I would like to thank Albert Lo for his heavy influence on the style and formulation of the calculus presented here. Thanks also to Hemant Ishwaran for his continuing help and encouragement.

DEPARTMENT OF INFORMATION
AND SYSTEMS MANAGEMENT
HONG KONG UNIVERSITY OF
SCIENCE AND TECHNOLOGY
CLEAR WATER BAY, KOWLOON
THE HONG KONG
E-MAIL: lancelot@ust.hk